\documentclass[letterpaper,12pt]{article}

\usepackage{amssymb,amsmath,amsthm,amsfonts}
\usepackage{xcolor}
\usepackage{hyperref}
\definecolor{linkblue}{rgb}{0.0,0.16,0.42}
\hypersetup{colorlinks=true,linkcolor=linkblue,citecolor=linkblue,urlcolor=linkblue}
\usepackage{graphicx}
\usepackage{tikz}
\usepackage{fullpage}
\usepackage{needspace}
\usepackage{xspace}
\usepackage{booktabs,tabularx}

\newcommand{\ord}{\mathrm{ord}}
\newcommand{\F}{\mathcal F}
\newcommand{\Tree}{\mathcal T}
\newcommand{\coef}[1]{\mathop{[\![#1]\!]}\xspace}
\newcommand{\im}{\mathop\mathrm{Im}\xspace}

\newcommand{\namedref}[2]{\hyperref[#2]{#1~\ref*{#2}}}
\newcommand{\namedeqref}[2]{\hyperref[#2]{#1~\textup{(\ref*{#2})}}}

\newtheorem{theorem}{Theorem}[section]
\newtheorem{lemma}[theorem]{Lemma}
\newtheorem{corollary}[theorem]{Corollary}
\newtheorem{proposition}[theorem]{Proposition}

\theoremstyle{definition}

\title{Counting Finite Maps through Functional Digraphs and Fibers}
\author{Max A. Alekseyev\thanks{The George Washington University, Washington, DC, United States. Email: \href{mailto:maxal@gwu.edu}{\texttt{maxal@gwu.edu}}.}}
\date{}

\begin{document}
\maketitle

\begin{abstract}
We count solutions of composition identities for maps on finite sets, first for three identities in which both maps vary and then for commutation with a prescribed arithmetic map. For $f\circ g=f\circ g\circ g$, drawing an arrow from each element to its image reduces the enumeration to rooted trees attached to directed cycles. For $f\circ f=f\circ g\circ f$, grouping elements with the same $f$-value leads instead to a partition-based enumeration. A classical rooted-forest enumeration leads in both problems to $\Tree(x)$, the exponential generating function for rooted labeled trees. For $f=g\circ f\circ g$, the values on the cycles of $g$ determine all remaining values. The exponential generating function $\mathcal D$ for these ordered pairs and the ordinary generating function $\mathcal A$ for endofunctions up to relabeling satisfy $\mathcal D(x)=\mathcal A(\Tree(x))$. Finally, we count functions on $\mathbb Z/n\mathbb Z$ satisfying $f(ax)=af(x)$, first for prime multipliers $a$ and then for a broader class of composite multipliers. In each problem, fixing one map exposes independent choices for the other through its functional digraph or its fibers.
\end{abstract}

\section{Introduction}

Unless the empty case is stated separately, $n$ is a positive integer, and we write $[n]:=\{1,2,\dots,n\}$.
Compositional identities for maps from $[n]$ to itself have rich combinatorial properties, and counting their solutions can be challenging.
Throughout, a map on a finite set $X$ means a self-map $X\to X$, also called an \emph{endofunction}.
Composition is read from right to left: $(fg)(x):=f(g(x))$, and $f^r$ denotes the $r$-fold iterate of $f$.
For example, pairs of commuting maps, satisfying $f\circ g=g\circ f$, were enumerated in~\cite{Holloway2013}.
A useful lesson is that there is no single best representation of a finite map: the form of the identity often suggests which structure should be exposed.

We first study
\begin{equation}\label{eq:fgfgg}
 f\circ g=f\circ g\circ g.
\end{equation}
For a map $g$, draw a directed arrow $x\to g(x)$ from every element to its image.
Repeatedly following arrows eventually reaches a directed cycle, while the other arrows form tree-like branches feeding into those cycles.
\namedref{Sections}{sec:graphs} and~\ref{sec:fgfgg} formalize this representation and reduce the enumeration for~\eqref{eq:fgfgg} to independent choices on sources and connected components.

We then consider
\begin{equation}\label{eq:fffgf}
 f\circ f=f\circ g\circ f.
\end{equation}
Here it is more effective to fix $f$ and group together elements having the same $f$-value.
Such a preimage $f^{-1}(y)$ is called a \emph{fiber} of $f$.
The nonempty fibers form a \emph{set partition}, i.e., a division of $[n]$ into disjoint nonempty blocks, and the condition on $g$ becomes a collection of independent choices governed by those blocks.
\namedref{Section}{sec:fibers} shows that a classical rooted-forest enumeration transforms the resulting weights into the same rooted-tree generating function that appears in the first problem.

Our third identity is
\begin{equation}\label{eq:gfg_intro}
 f=g\circ f\circ g.
\end{equation}
A simple change of variable shows that the values on points lying on directed cycles under repeated application of $g$ determine the values everywhere else.
When both maps are permutations, the count is $n!p(n)$, where $p(n)$ is the number of \emph{integer partitions} of $n$, i.e., unordered sums of positive integers with total $n$.
For arbitrary maps, \namedref{Section}{sec:gfg} relates the corresponding count to endofunctions considered up to relabeling and then determines how the remaining points can be attached.

Finally, for a positive integer $a$ we study
\[
 f(ax)=a f(x)
\]
on the residue classes modulo $n$.
Here $a$ and $n$ are prescribed, and only $f$ varies: the equation asks for maps that commute with multiplication by $a$.
\namedref{Section}{sec:px} first treats prime multipliers, for which the points outside the directed cycles have a regular branching structure. For a squarefree composite multiplier, the same structure survives when the primes common to $a$ and $n$ occur in $n$ with equal exponents.
The first three problems also begin by fixing one map and counting its compatible partners.
When both maps vary, we then sum these counts; in the arithmetic application, the prescribed map itself supplies additional structure.
Thus the four problems illustrate a common strategy: choose a representation in which the equation separates into independent choices.

Numerical values for the first three counting problems are collected in the On-Line Encyclopedia of Integer Sequences (OEIS) as A239785, A239772, and A239753; the case $a=2$ of the fourth problem gives A117987~\cite{OEIS}.
These entries provide numerical data for comparison with the formulas derived below.
When both maps are permutations, cancellation in either of the first two identities forces $g$ to be the identity map, giving $n!$ pairs in each case.
For the third problem, the permutation subcase is classical; \namedref{Section}{sec:gfg} relates it to pairs of permutations with equal squares before treating arbitrary maps.

\medskip
\noindent\textit{Surrounding identities.}
Short identities provide a useful testing ground: changing the order of composition can change which information about a map matters.
Call an identity of \emph{type $(\ell,r)$} one whose two sides are compositions of $\ell$ and $r$ maps chosen from $f,g$, with $\ell\leq r$.
Interchanging the names $f,g$ normalizes the shorter side of a type $(1,3)$ identity to $f$, leaving the eight forms in the table.
Reversing the order of composition is not included in this normalization.
A composition is called \emph{mixed} if it uses both $f$ and $g$.
The five mixed cases other than $f=gfg$ have explicit image and cycle enumerations in Einstein's mathematical notebooks~\cite{EinsteinNotebooks}.
For $f=f^3$, the restriction of $f$ to its image is an involution, meaning a permutation whose square is the identity, and $g$ is arbitrary; for $f=g^3$, every $g$ determines a unique $f$.
Together with \namedref{Section}{sec:gfg}, these descriptions give explicit enumerations for all eight forms.

\begin{table}[htbp]
\centering
\small
\renewcommand{\arraystretch}{1.12}
\begin{tabularx}{\textwidth}{@{}l l >{\raggedright\arraybackslash}X@{}}
\toprule
Identity & Numerical data & Counting description\\
\midrule
$f=f^3$ & $n^n\cdot\mathrm{A060905}(n)$ & An involution on the image; $g$ is free.\\
$f=f^2g$ & A239760 & Weighted cycles with leaves.\\
$f=fgf$ & A239768 & An injective choice of one representative from each fiber.\\
$f=fg^2$ & A239777 & Constant values on the components of $g^2$.\\
$f=gf^2$ & A235328 & $f$ permutes its image; $g$ is its inverse there.\\
$f=gfg$ & A239753 & Cycle assignments and rooted forests; \namedref{Section}{sec:gfg}.\\
$f=g^2f$ & A239771 & $g^2$ fixes the image of $f$ pointwise.\\
$f=g^3$ & A000312 & The count is $n^n$.\\
\bottomrule
\end{tabularx}
\caption{The eight normalized type $(1,3)$ identities. The middle column identifies numerical data in the OEIS~\cite{OEIS}; A060905 counts the admissible maps $f$ alone.}
\label{tab:type13}
\end{table}

For type $(2,3)$, the normalized shorter side is $f^2$ or $fg$, giving sixteen forms.
Fifteen involve both maps, while $f^2=f^3$ is a unary condition with the other map unrestricted.
Fourteen have at least one mixed composition; the two exceptions are $f^2=f^3$ and the equal-power equation $f^2=g^3$.
A follow-up study~\cite{AlekseyevType23} combines these viewpoints with additional methods to count solutions of all sixteen identities.
The two examples treated here exhibit the same rooted-tree generating function, which also enters the type $(1,3)$ count in \namedref{Section}{sec:gfg}.
The modular equation has a separate role: it is a type $(2,2)$ commutation relation with multiplication by $a$ prescribed, and provides an arithmetic application of the same functional-digraph methods.

\section{Main results}\label{sec:main_results}
We state the four principal enumeration theorems, together with the generating-function notation that accompanies them.
For a formal power or Laurent series $F$, $\coef{x^n}F(x)$ denotes the coefficient of $x^n$ in $F(x)$; the same notation is used for multivariate coefficient extraction, with $\coef{z^0}$ denoting the constant term in $z$.
For counts $c_r$, the \emph{ordinary generating function} is $\sum_{r\geq0}c_r x^r$, while the \emph{exponential generating function} (EGF) is $\sum_{r\geq0}c_r x^r/r!$.
Define
\begin{equation}\label{eq:Tree_intro_def}
 \Tree(x):=\sum_{r\geq1} r^{r-1}\frac{x^r}{r!}=-W(-x),
\end{equation}
where $W$ denotes the Lambert $W$ function.  The second equality follows directly from its Taylor expansion at the origin,
\[
 W(x)=\sum_{r\geq1}(-r)^{r-1}\frac{x^r}{r!};
\]
see~\cite{Corless1996}.  \namedref{Section}{sec:graphs} gives $\Tree(x)$ its combinatorial interpretation as the exponential generating function for rooted labeled trees and derives its basic functional equation.
We also write $\mathbb Z_n:=\mathbb Z/n\mathbb Z$.

\begin{theorem}\label{th:number_fg}
Let $A_n$ be the number of ordered pairs of maps $(f,g)$ on $[n]$ satisfying~\eqref{eq:fgfgg}.  Then
\begin{equation}\label{eq:An}
 A_n=n!\,\coef{x^n}\left(1-\Tree\!\left(xe^{(n-1)x}\right)\right)^{-n}.
\end{equation}
\end{theorem}
We prove \namedref{Theorem}{th:number_fg} in \namedref{Section}{sec:fgfgg}.

\begin{theorem}\label{th:number_fffgf}
Let $C_n$ be the number of ordered pairs of maps $(f,g)$ on $[n]$ satisfying~\eqref{eq:fffgf}.  Then
\begin{equation}\label{eq:Cn_constant_term}
 C_n=n!\,\coef{z^0 y^n}
 \left(
  1-\sum_{i=1}^n
  \frac{y^i(n+iz)^i}{i!\,z}
 \right)^{-1}.
\end{equation}
The inverse is expanded in $\mathbb Q[z,z^{-1}][[y]]$, and the displayed operator extracts the coefficient of $y^n z^0$.
\end{theorem}
We prove \namedref{Theorem}{th:number_fffgf} in \namedref{Section}{sec:fibers}, first counting by fiber sizes and then collecting those counts into~\eqref{eq:Cn_constant_term}. A rooted-forest identity expresses the resulting fiber series in the same tree function $\Tree(x)$. \namedref{Appendix}{sec:Cn_univariate} gives an alternative one-variable formula.

For the third identity, unlabeled endofunctions enter naturally.
Let $a_k$ be the number of endofunctions on a $k$-element set up to relabeling, with $a_0:=1$, and define their ordinary generating function
\[
 \mathcal A(u):=\sum_{k\geq0}a_k u^k.
\]
Equivalently, $a_k$ is the number of unlabeled functional digraphs on $k$ vertices.
Let $\mathcal R(u)$ denote the ordinary generating function for unlabeled rooted trees, so the coefficient of $u^k$ counts such trees on $k$ vertices.
The classical cycle-and-tree decomposition of unlabeled functional digraphs gives the product below; see Harary~\cite{Harary1959} and~\cite[Note VII.25, p.~480]{FlajoletSedgewick2009}:
\begin{equation}\label{eq:A_R_product}
 \mathcal A(u)=\prod_{j\geq1}\frac{1}{1-\mathcal R(u^j)}.
\end{equation}
The series $\mathcal R(u)$ satisfies the functional identity
\begin{equation}\label{eq:Polya_rooted_tree}
 \mathcal R(u)
 =u\exp\!\left(\sum_{j\geq1}\frac{\mathcal R(u^j)}{j}\right).
\end{equation}

\Needspace{9\baselineskip}
\begin{theorem}\label{th:number_gfg}
Let $D_n$ be the number of ordered pairs of maps $(f,g)$ on $[n]$ satisfying~\eqref{eq:gfg_intro}.
Then $D_0=1$ and, for $n\geq1$,
\begin{equation}\label{eq:Dn_gfg}
 D_n=\sum_{k=1}^n\binom{n}{k}\,k n^{n-k-1}k!a_k.
\end{equation}
Equivalently, if $\mathcal D(x):=\sum_{n\geq0}D_nx^n/n!$, then
\begin{equation}\label{eq:D_composition}
 \mathcal D(x)=\mathcal A(\Tree(x)).
\end{equation}
\end{theorem}
We prove \namedref{Theorem}{th:number_gfg} in \namedref{Section}{sec:gfg}.  The permutation subcase, also treated there, is $n!p(n)$, where $p(n)$ is the number of integer partitions of $n$.

For the fourth result, let $p$ be prime and write $n=p^tm$ with $p\nmid m$.
Let $\mu$ be the arithmetic M\"obius function, defined by $\mu(q)=0$ when $q$ is divisible by the square of a prime and $\mu(q)=(-1)^s$ when $q$ is a product of $s$ distinct primes.
If $m>1$, let $r$ be the least positive integer such that $p^r\equiv1\pmod m$, and define
\begin{equation}\label{eq:Nd}
 N_d:=\sum_{q\mid d}\mu(d/q)\gcd(m,p^q-1).
\end{equation}

\begin{theorem}\label{th:number_px}
If $m>1$, the number $B_p(n)$ of maps $f$ on $\mathbb Z_n$ satisfying $f(px)=pf(x)$ is
\begin{equation}\label{eq:Bpn}
 B_p(n)=p^{n-m}\prod_{d\mid r}\gcd(m,p^d-1)^{N_d/d}.
\end{equation}
If $m=1$, then $B_p(n)=p^{n-1}$.
\end{theorem}
We prove \namedref{Theorem}{th:number_px} in \namedref{Section}{sec:px}, where \namedref{Theorem}{th:number_balanced} gives a composite-multiplier extension.

\section{Functional digraphs and rooted forests}\label{sec:graphs}

The \emph{functional digraph} of a map $h:[n]\longrightarrow[n]$, i.e., the directed graph recording the action of $h$, is the labeled directed graph $G(h)$ with vertex set $[n]$ and arc set
\[
 \{(i,h(i)):i\in[n]\}.
\]
Thus every vertex of $G(h)$ has outdegree $1$.
A vertex of indegree $0$ will be called a \emph{source}.
A loop corresponds to a fixed point of $h$.
Throughout the paper, by a connected component of a functional digraph we mean a connected component of its underlying undirected graph.

Starting at any vertex and repeatedly following its unique outgoing arc, we eventually reach a directed cycle, possibly a loop.
Consequently, every connected component of $G(h)$ contains exactly one directed cycle.
A vertex lying on a directed cycle will be called \emph{cyclic}; its \emph{period} is the length of that cycle.
The union $C(h)$ of these cycles will be called the \emph{core} of $h$; the restriction of $h$ to the vertices of $C(h)$ is a permutation.
When $h$ itself is a permutation, every vertex belongs to $C(h)$, and these directed cycles are exactly the cycles in the usual disjoint-cycle decomposition of the permutation.
This standard decomposition of a finite mapping into directed cycles with rooted labeled trees attached to their vertices is discussed, for example, in~\cite[Section II.5]{FlajoletSedgewick2009}.

A \emph{rooted forest} is a collection of rooted trees.  Deleting the cycle arcs from $G(h)$ produces a \emph{planted forest}, i.e., a forest of rooted directed trees whose arcs point toward their roots; we denote it by $T(h)$.
Each tree contains exactly one vertex of $C(h)$, which serves as its root.
A \emph{leaf}, i.e., a vertex with no children in its rooted tree, is a source of $G(h)$ whenever it is not a root.
A singleton tree is the only exception to this leaf/source correspondence: its root is a leaf of the tree but is not a source of $G(h)$, because it has an incoming cycle arc in the original functional digraph.

\namedref{Figure}{fig:functional-digraph} illustrates this decomposition.
For the displayed map, the sources are $5$ and $8$; after the cycle arcs are removed, the core vertices become the roots of the planted forest.

\begin{figure}[ht]
\centering
\begin{tikzpicture}[>=stealth, every node/.style={circle,draw,minimum size=7mm,inner sep=0pt}]
  \node (v1) at (0,0) {1};
  \node (v2) at (1.6,0) {2};
  \node (v3) at (0.8,-1.4) {3};
  \node (v4) at (0.8,-2.8) {4};
  \node (v5) at (0.8,-4.2) {5};
  \node (v6) at (4.2,0) {6};
  \node (v7) at (5.8,0) {7};
  \node (v8) at (7.2,0) {8};

  \draw[->,thick] (v1) -- (v2);
  \draw[->,thick] (v2) -- (v3);
  \draw[->,thick] (v3) -- (v1);
  \draw[->] (v4) -- (v3);
  \draw[->] (v5) -- (v4);
  \draw[->,thick,bend left=18] (v6) to (v7);
  \draw[->,thick,bend left=18] (v7) to (v6);
  \draw[->] (v8) -- (v7);
\end{tikzpicture}
\caption{A functional digraph. Thick arcs form the core $C(h)$; deleting them leaves the planted forest $T(h)$.}
\label{fig:functional-digraph}
\end{figure}
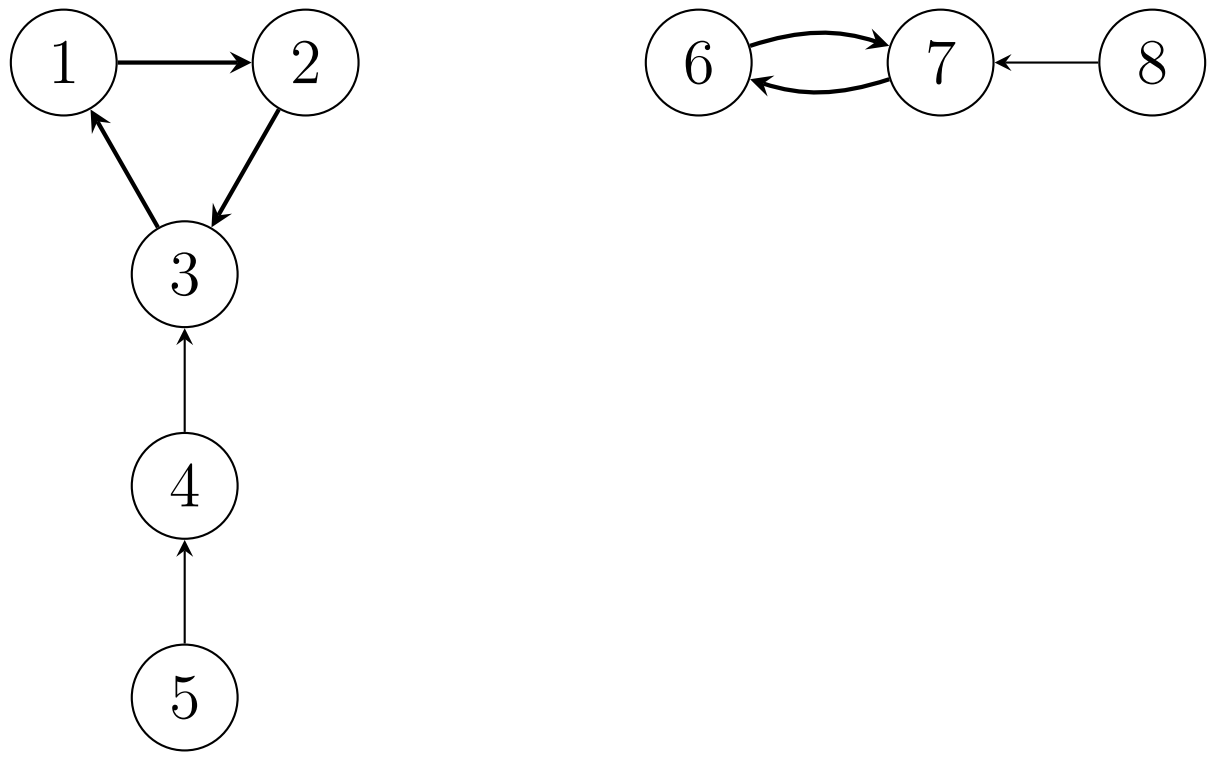

We will use one further piece of general machinery throughout the paper: the generating function for rooted labeled trees, refined when necessary by leaf weights.

Let $t_r$ be the number of rooted labeled directed trees on $r$ vertices, with every arc directed toward the root.
Cayley's rooted-tree formula gives $t_r=r^{r-1}$; see, for example,~\cite[Section II.5]{FlajoletSedgewick2009}.
Hence the EGF of these trees is exactly the series $\Tree(x)$ defined in~\eqref{eq:Tree_intro_def}.
A rooted tree consists recursively of a root together with a set of rooted trees attached to it.
The labeled exponential formula~\cite[Section 5.1]{Stanley1999v2} translates this recursive decomposition into
\begin{equation}\label{eq:tree_function}
 \Tree(x)=x e^{\Tree(x)}.
\end{equation}
This series is the classical \emph{tree function}.

The functional equation also gives another explanation of the Lambert $W$ representation in~\eqref{eq:Tree_intro_def}.  Since $W$ is characterized near $0$ by
$W(z)e^{W(z)}=z$ and $W(0)=0$, we have $-\Tree(x)e^{-\Tree(x)}=-x$, and therefore
\begin{equation}\label{eq:tree_W}
 \Tree(x)=-W(-x).
\end{equation}
We keep $\Tree(x)$ as the primary notation because it records the combinatorial origin of the series, while~\eqref{eq:tree_W} supplies the equivalent Lambert $W$ form.

For our application we must also record leaves.
Let $t_r(y)$ be the polynomial in which the coefficient of $y^j$ is the number of rooted directed trees on $r$ labeled vertices with $j$ leaves, counting the root as a leaf when the tree consists of a single vertex, and let
\[
 \F(x,y):=\sum_{r=1}^{\infty} t_r(y)\frac{x^r}{r!}.
\]

\begin{theorem}\label{th:gfFxy}
The bivariate exponential generating function $\F(x,y)$ is
\[
 \F(x,y)=\Tree\!\left(xe^{x(y-1)}\right)+x(y-1)
        =-W\!\left(-xe^{x(y-1)}\right)+x(y-1).
\]
\end{theorem}

\begin{proof}
The same recursive decomposition that gives~\eqref{eq:tree_function} applies to the leaf-weighted series $\F(x,y)$, with one correction.
In~\eqref{eq:tree_function}, the term $x$ corresponding to a root with no attached trees represents the singleton rooted tree.
For $\F(x,y)$ this tree has one leaf, so its contribution must be $xy$ instead of $x$.
All nonsingleton rooted trees are counted by the same recursive construction, since their roots have at least one child and hence are not leaves.
Thus the required correction is to add $x(y-1)$ to the right-hand side of~\eqref{eq:tree_function}, with $\Tree(x)$ replaced by $\F(x,y)$.
Hence
\[
 \F(x,y)=x e^{\F(x,y)}+x(y-1),
\]
which we rewrite as
\[
 \F(x,y)-x(y-1)  = xe^{x(y-1)} e^{\F(x,y)-x(y-1)}.
\]
Then, by~\eqref{eq:tree_function}, we have
\[
 \F(x,y)-x(y-1)
 =\Tree\!\left(xe^{x(y-1)}\right)
 =-W\!\left(-xe^{x(y-1)}\right),
\]
which proves the formula for $\F(x,y)$.
\end{proof}

For a fixed positive integer $N$, putting $y=N$ gives the convenient one-variable abbreviation
\begin{equation}\label{eq:HN}
 H_N(x):=\F(x,N)-x(N-1)
       =\Tree\!\left(xe^{(N-1)x}\right)
       =-W\!\left(-xe^{(N-1)x}\right).
\end{equation}

Fix positive integers $N,r,$ and $k$.
For a planted forest $T$ on a prescribed $r$-element label set with exactly $k$ rooted trees, let $s(T)$ denote the number of nonroot leaves and give $T$ the weight $N^{s(T)}$.

\begin{theorem}\label{th:gf_forest}
The total weight of the planted forests specified above is
\[
 \sum_T N^{s(T)}=r!\,\coef{x^r}\frac{H_N(x)^k}{k!}.
\]
\end{theorem}

\begin{proof}
By~\eqref{eq:HN}, $H_N(x)$ is the exponential generating function for one rooted labeled tree in which each nonroot leaf has weight $N$.
A forest with exactly $k$ rooted components is an unordered set of $k$ such trees.
The labeled exponential formula therefore gives $H_N(x)^k/k!$.
Extracting the coefficient of $x^r/r!$ proves the claim.
\end{proof}

\section{The identity \texorpdfstring{$f\circ g=f\circ g\circ g$}{f o g = f o g o g} on \texorpdfstring{$[n]$}{[n]}}\label{sec:fgfgg}

The following theorem is the key reduction.
For a map $g$ on $[n]$, let $s(g)$ and $c(g)$ denote the number of sources and connected components of $G(g)$, respectively.

\begin{theorem}\label{th:fixed_g}
For a fixed map $g$ on $[n]$, exactly
\[
 n^{s(g)+c(g)}
\]
maps $f$ on $[n]$ satisfy $f\circ g=f\circ g\circ g$.
\end{theorem}

\begin{proof}
A source of $G(g)$ is an element of $[n]$ that is not in the image of $g$.
The value of $f$ on such an element never appears in~\eqref{eq:fgfgg}, because there $f$ is always applied to a value of $g$.
Hence $f$ may take arbitrary and independent values on the $s(g)$ sources, contributing a factor $n^{s(g)}$.

Now let $v$ be a non-source vertex.
Choose an incoming arc $(u,v)$, so that $g(u)=v$.
Then~\eqref{eq:fgfgg} gives
\[
 f(v)=f(g(u))=f(g^2(u))=f(g(v)).
\]
Thus $f$ takes the same value at a non-source vertex and at its successor under $g$.
Following successors shows that $f$ is constant around the directed cycle of each component.
Every non-source vertex eventually reaches that cycle, and applying the same equality along its successive iterates shows that it has the same $f$-value as the cycle vertices.
Therefore $f$ is constant on all non-source vertices of each connected component.
Conversely, every such assignment satisfies the identity, because $g(x)$ and $g^2(x)$ are non-source vertices in the same component.

Every component contains at least one non-source vertex, namely every vertex of its directed cycle.
Hence $f$ is completely determined by its independent values on the sources and by one additional value for the non-source vertices in each of the $c(g)$ components.
The total number of choices is $n^{s(g)+c(g)}$.
\end{proof}

We can now prove the first main theorem.

\begin{proof}[Proof of \namedref{Theorem}{th:number_fg}]
Let $k$ be the number of trees in $T(g)$, equivalently the number of vertices in the core $C(g)$.
By \namedref{Theorem}{th:gf_forest}, with $N=r=n$, the weighted contribution of all forests with $k$ trees, where the weight records the factor $n^{s(g)}$ from \namedref{Theorem}{th:fixed_g}, is
\[
 n!\,\coef{x^n}\frac{H_n(x)^k}{k!},
 \qquad
 H_n(x)=\Tree\!\left(xe^{(n-1)x}\right).
\]

Once the $k$ rooted trees are fixed, the cycle arcs of $g$ are specified by a permutation of their roots.
If this permutation has $t$ cycles in its disjoint-cycle decomposition, then $G(g)$ has $t$ connected components.
There are $\genfrac{[}{]}{0pt}{}{k}{t}$ such permutations, where $\genfrac{[}{]}{0pt}{}{k}{t}$ is the unsigned Stirling number of the first kind.
Hence the total contribution of the factor $n^{c(g)}$ over all permutations of the $k$ roots is
\[
 \sum_{t=1}^k \genfrac{[}{]}{0pt}{}{k}{t}n^t
 =n(n+1)\cdots(n+k-1)
 =\frac{(n+k-1)!}{(n-1)!}.
\]
Therefore
\begin{align*}
 A_n
 &=n!\,\coef{x^n}\sum_{k=1}^n
   \frac{(n+k-1)!}{(n-1)!}\frac{H_n(x)^k}{k!}\\
 &=n!\,\coef{x^n}\sum_{k=1}^n
   \binom{n+k-1}{k}H_n(x)^k.
\end{align*}
Since $H_n(x)$ has zero constant term, terms with $k>n$ cannot contribute to $\coef{x^n}$.
We may therefore extend the sum to all $k\geq0$ and use
$\sum_{k\geq0}\binom{n+k-1}{k}u^k=(1-u)^{-n}$, obtaining~\eqref{eq:An}.
\end{proof}

The first few values are
\[
\begin{aligned}
 A_1&=1, & A_2&=14, & A_3&=411,\\
 A_4&=19912, & A_5&=1412745, & A_6&=136537056.
\end{aligned}
\]
These values form OEIS sequence A239785~\cite{OEIS}, providing a numerical illustration of~\eqref{eq:fgfgg}.

\section{The identity \texorpdfstring{$f\circ f=f\circ g\circ f$}{f o f = f o g o f}: a fiber viewpoint}\label{sec:fibers}

The preceding problem was naturally organized by fixing $g$ and examining its functional digraph.
For the identity~\eqref{eq:fffgf}, it is more effective to reverse the point of view and fix $f$.
The reason is that the inner occurrence of $f$ restricts the identity to the range of $f$, while the outer occurrence of $f$ turns the condition on $g$ into membership in a fiber of $f$.

Let $k:=|\im f|$.
For $y\in \im f$, the \emph{fiber} of $f$ over $y$, i.e., the preimage $f^{-1}(y)$, is nonempty.
These fibers form a set partition of $[n]$ into $k$ nonempty blocks.

\begin{lemma}\label{lem:fixed_f}
For a fixed map $f:[n]\to[n]$ with $|\im f|=k$, the number of maps $g:[n]\to[n]$ satisfying
$f\circ f=f\circ g\circ f$ is
\begin{equation}\label{eq:fixed_f_count}
 n^{n-k}\prod_{y\in \im f}\bigl|f^{-1}(f(y))\bigr|.
\end{equation}
\end{lemma}

\begin{proof}
As $x$ ranges over $[n]$, the value $y=f(x)$ ranges over $\im f$.
Hence~\eqref{eq:fffgf} is equivalent to
\[
 f(y)=f(g(y))\qquad (y\in \im f).
\]
For a fixed $y\in \im f$, this condition says precisely that
\[
 g(y)\in f^{-1}(f(y)),
\]
so there are $|f^{-1}(f(y))|$ possible values of $g(y)$.
These choices are independent for the $k$ elements of $\im f$.
On the remaining $n-k$ elements of $[n]\setminus \im f$, the map $g$ is unrestricted, contributing the factor $n^{n-k}$.
Multiplying the independent choices gives~\eqref{eq:fixed_f_count}.
\end{proof}

For example, \namedref{Figure}{fig:fibers} shows a map with fibers
$f^{-1}(2)=\{1,3\}$, $f^{-1}(4)=\{2,4,6\}$, and $f^{-1}(6)=\{5\}$.
Its image is $\{2,4,6\}$, and each of these three points belongs to the middle fiber.
Thus each of $g(2),g(4),g(6)$ has three choices, while $g(1),g(3),g(5)$ are unrestricted, giving $3^3 6^3=5832$ choices in all.
The condition constrains $g$ at image points; it does not require $g$ to preserve whole fibers.

\begin{figure}[htbp]
\centering
\begin{tikzpicture}[>=stealth,
  point/.style={circle,draw,fill=white,minimum size=6.5mm,inner sep=0pt}]
  \draw (0,0) ellipse (1.25 and 1.05);
  \draw[fill=black!5] (3.5,0) ellipse (1.45 and 1.05);
  \draw (7,0) ellipse (1.05 and 1.05);
  \node at (0,1.4) {$f^{-1}(2)$};
  \node at (3.5,1.4) {$f^{-1}(4)$};
  \node at (7,1.4) {$f^{-1}(6)$};
  \node[point] (p1) at (-0.5,0) {$1$};
  \node[point] (p3) at (0.5,0) {$3$};
  \node[point] (p2) at (2.9,0.25) {$\mathbf{2}$};
  \node[point] (p4) at (4.05,0.25) {$\mathbf{4}$};
  \node[point] (p6) at (3.5,-0.55) {$\mathbf{6}$};
  \node[point] (p5) at (7,0) {$5$};
  \draw[->,dashed] (p2) -- (p4);
  \draw[->,dashed] (p6) -- (p4);
  \draw[->,dashed] (p4) to[loop right,min distance=5mm] (p4);
\end{tikzpicture}
\caption{Fibers of $f$ on $[6]$. Bold labels form $\im f$. The dashed arrows show the permitted choice $g(2)=g(4)=g(6)=4$; the other values of $g$ are free.}
\label{fig:fibers}
\end{figure}
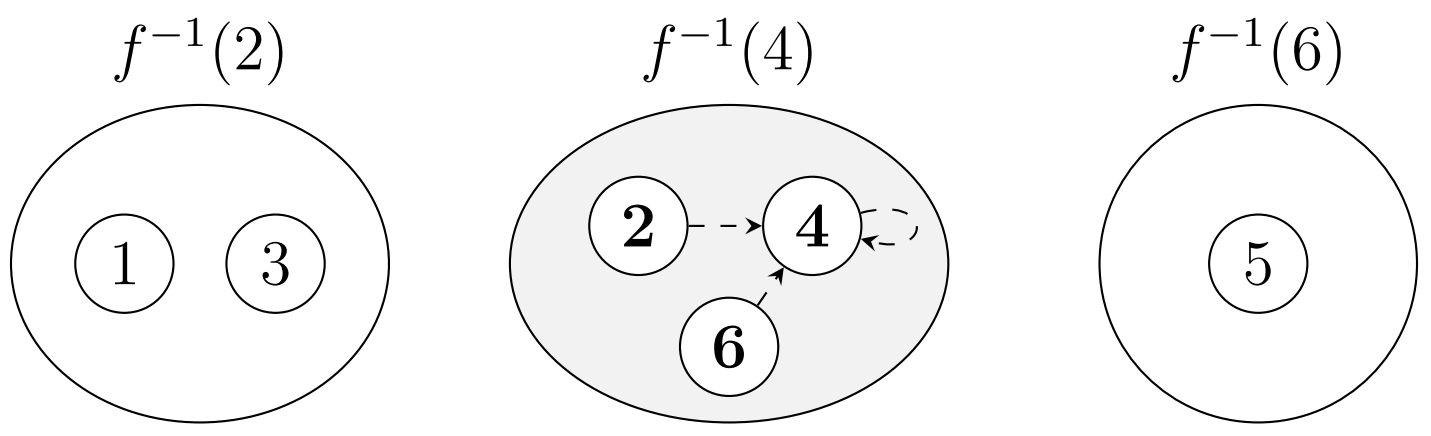

\subsection{Enumeration by fiber profile}

For $i=1,\dots,n$, let
\[
 m_i:=\#\{y\in \im f:|f^{-1}(y)|=i\}.
\]
Thus $m_i$ is the number of nonempty fibers of size $i$, and
\begin{equation}\label{eq:fiber_profile}
 \sum_{i=1}^n i m_i=n,
 \qquad
 \sum_{i=1}^n m_i=k.
\end{equation}
The fibers form a set partition of $[n]$, and we refer to $(m_1,\dots,m_n)$ as its \emph{profile}.

A set partition with a given profile $(m_1,\dots,m_n)$ can be chosen in
\begin{equation}\label{eq:number_partitions_profile}
 \frac{n!}{\displaystyle\prod_{i=1}^n m_i!(i!)^{m_i}}
\end{equation}
ways.
Fix such a set partition, and suppose that $k_i$ elements of the range $\im f$ lie in blocks of size $i$.
Then
\[
 0\leq k_i\leq i m_i,
 \qquad
 \sum_{i=1}^n k_i=k.
\]
There are
\[
 \prod_{i=1}^n\binom{i m_i}{k_i}
\]
ways to choose these range elements.
A point in a block of size $i$ contributes $1$ if it is omitted from the range, or $iz$ if it is selected: $z$ marks the selected image value, and $i$ counts its permitted images under $g$ by \namedref{Lemma}{lem:fixed_f}.
This weight depends on the block containing the point, not on the later assignment of image values to blocks.
Thus the $im_i$ points in blocks of size $i$ contribute $(1+iz)^{im_i}$, and the weighted number of choices for the $k$-element range is
\begin{equation}\label{eq:range_weight}
 \coef{z^k}\prod_{i=1}^n(1+iz)^{i m_i}.
\end{equation}

Once the set partition and its $k$-element range have been fixed, there are $k!$ ways to assign the $k$ blocks bijectively to the $k$ range values.
Each such bijection determines a unique map $f$ with the prescribed fibers and range.
Combining these observations gives the fiber-profile form of the second main enumeration.

\begin{theorem}\label{th:number_fffgf_profile}
Let $C_n$ be the number of ordered pairs of maps $(f,g)$ on $[n]$ satisfying
$f\circ f=f\circ g\circ f$.
For each profile $(m_1,\dots,m_n)$ in the sum below, put $k:=\sum_{i=1}^n m_i$.
Then
\begin{equation}\label{eq:Cn_profile}
 C_n=
 n!\!\sum_{\substack{m_1,\dots,m_n\geq0\\
                       \sum_{i=1}^n i m_i=n}}
 \frac{k!\,n^{n-k}}
      {\displaystyle\prod_{i=1}^n m_i!(i!)^{m_i}}
 \coef{z^k}\prod_{i=1}^n(1+iz)^{i m_i}.
\end{equation}
\end{theorem}

\begin{proof}
For a fixed fiber profile, the factors in~\eqref{eq:Cn_profile} have the following meanings.
\namedeqref{Equation}{eq:number_partitions_profile} chooses the set partition into fibers; the coefficient in~\eqref{eq:range_weight} chooses the range and simultaneously records the product in \namedref{Lemma}{lem:fixed_f}; the factor $k!$ assigns the fibers bijectively to the range values; and $n^{n-k}$ accounts for the unrestricted values of $g$ outside the range of $f$.
Summing over all profiles satisfying~\eqref{eq:fiber_profile} counts every pair $(f,g)$ exactly once.
\end{proof}

Only the sizes of the blocks containing the selected image values enter their weights.
We can therefore collect fibers of each size into a generating function and sum over all profiles at once.

\subsection{A constant-term form and the tree function}

The variable $y$ will record the total number of points, while $z$ records selected image values. Collecting the weights by fiber size turns the profile sum into a coefficient extraction.
For fixed $k$, the multinomial theorem gives
\[
 \sum_{\substack{m_1,m_2,\dots\geq0\\\sum_i m_i=k}}
 \frac{k!}{\prod_i m_i!}
 \prod_{i\geq1}
 \left(\frac{y^i(1+iz)^i}{i!}\right)^{m_i}
 =
 \left(\sum_{i\geq1}\frac{y^i(1+iz)^i}{i!}\right)^k.
\]
Using this in~\eqref{eq:Cn_profile} yields
\begin{equation}\label{eq:Cn_diagonal}
 C_n=
 n!\sum_{k=1}^n n^{n-k}
 \coef{y^n z^k}
 \left(\sum_{i\geq1}\frac{y^i(1+iz)^i}{i!}\right)^k.
\end{equation}
Now make the simultaneous rescaling $y\mapsto ny$ and $z\mapsto z/n$.
For every $k$,
\[
 n^{n-k}
 \coef{y^n z^k}
 \left(\sum_{i\geq1}\frac{y^i(1+iz)^i}{i!}\right)^k
 =
 \coef{y^n z^k}
 \left(\sum_{i\geq1}\frac{y^i(n+iz)^i}{i!}\right)^k.
\]
Define
\begin{equation}\label{eq:Qn}
 Q_n(y,z):=\sum_{i\geq1}\frac{y^i(n+iz)^i}{i!}.
\end{equation}
This rescaled fiber series already has a closed form in terms of the tree function:
\begin{equation}\label{eq:Qn_tree}
 1+Q_n(y,z)
 =\frac{\exp\!\bigl(n\Tree(yz)/z\bigr)}{1-\Tree(yz)}
 =\frac{\exp\!\bigl(-nW(-yz)/z\bigr)}{1+W(-yz)}.
\end{equation}
Thus the same rooted-tree series that appeared in the functional-digraph enumeration is already present in the rescaled fiber weights.

To verify~\eqref{eq:Qn_tree}, we use a classical rooted-forest identity related to the Abel polynomials; see~\cite[Exercise 5.37(b)]{Stanley1999v2}.
For a positive integer $a$, the exponential generating function $e^{a\Tree(u)}$ counts rooted labeled forests in which each component receives one of $a$ colors.
Equivalently, adjoin $a$ distinguished vertices, one for each color, and attach the root of every component to the distinguished vertex of its color.
The rooted-forest form of Cayley's formula then gives $a(a+i)^{i-1}$ such forests on $i$ ordinary labeled vertices, so
\begin{equation}\label{eq:abel_forest}
 e^{a\Tree(u)}
 =\sum_{i\geq0}a(a+i)^{i-1}\frac{u^i}{i!}.
\end{equation}
Since each coefficient is a polynomial in $a$, this identity also holds formally in $a$.
Differentiating~\eqref{eq:tree_function} gives
\[
 u\Tree'(u)=\frac{\Tree(u)}{1-\Tree(u)}.
\]
Consequently,
\begin{align}
 \sum_{i\geq0}(a+i)^i\frac{u^i}{i!}
 &=e^{a\Tree(u)}+\frac{u}{a}\frac{d}{du}e^{a\Tree(u)} \notag\\
 &=\frac{e^{a\Tree(u)}}{1-\Tree(u)}.\label{eq:abel_tree_identity}
\end{align}
Putting $a:=n/z$ and $u:=yz$ gives~\eqref{eq:Qn_tree}, since
$u^i(a+i)^i=y^i(n+iz)^i$.
Although $a=n/z$ is displayed with a negative power of $z$, the quotient $\Tree(yz)/z$ lies in $y\mathbb Q[[y,z]]$, so the right-hand side of~\eqref{eq:Qn_tree} is an ordinary formal power series in $y$ and $z$.

In $Q_n(y,z)^k/z^k$, the exponent of $z$ is the number of selected image values minus the number of fibers. Taking the constant term makes these numbers equal:
\[
 \coef{y^n z^k}Q_n(y,z)^k
 =\coef{z^0 y^n}\left(\frac{Q_n(y,z)}{z}\right)^k.
\]
Since $Q_n(y,z)$ has zero constant term in $y$, the sum over $k$ is geometric; the $k=0$ term does not contribute to $\coef{y^n}$ for $n\geq1$.
Thus
\begin{equation}\label{eq:Cn_constant_Q}
 C_n=n!\,\coef{z^0 y^n}
 \left(1-\frac{Q_n(y,z)}{z}\right)^{-1}.
\end{equation}
The inverse is expanded in $\mathbb Q[z,z^{-1}][[y]]$.
Terms with $i>n$ in~\eqref{eq:Qn} cannot affect $\coef{y^n}$, so~\eqref{eq:Cn_constant_Q} is exactly~\eqref{eq:Cn_constant_term} and proves \namedref{Theorem}{th:number_fffgf}.
Using~\eqref{eq:Qn_tree}, the same formula can equivalently be written as
\begin{equation}\label{eq:Cn_tree_constant}
 C_n=n!\,\coef{z^0 y^n}
 \left[
  1-\frac1z
  \left(
   \frac{\exp\!\bigl(n\Tree(yz)/z\bigr)}{1-\Tree(yz)}-1
  \right)
 \right]^{-1}.
\end{equation}

The constant term can also be eliminated; \namedref{Appendix}{sec:Cn_univariate} gives the resulting one-variable formula and its proof.

The first values are
\[
 C_1=1,\quad C_2=10,\quad C_3=231,\quad C_4=9688,
 \quad C_5=603445,\quad C_6=52284816,
\]
followed by $C_7=5951141035$ and $C_8=856275088768$.
Together with the conventional value $C_0=1$, these form OEIS sequence A239772~\cite{OEIS}, corresponding to the identity
$f(f(x))=f(g(f(x)))$.

The first two compositional identities therefore suggest complementary structural descriptions.
For~\eqref{eq:fgfgg}, the functional digraph of $g$ exposes the condition; for~\eqref{eq:fffgf}, the fiber partition of $f$ does so.
Our next identity returns to functional digraphs, but with an especially simple feature: only the cycle core carries any freedom.

\section{The identity \texorpdfstring{$f=g\circ f\circ g$}{f = g o f o g}: cycles, relabeling, and rooted trees}\label{sec:gfg}

Consider
\begin{equation}\label{eq:gfg}
 f=g\circ f\circ g
\end{equation}
on $[n]$.
Fix a map $g:[n]\to[n]$ and write
\[
 \mathcal S_g:=\{f:[n]\to[n]:f=g\circ f\circ g\}.
\]
On this solution set, define left and right composition by $g$ as
\[
 L_g(f):=g\circ f,
 \qquad
 R_g(f):=f\circ g.
\]
The first observation is that these operations permute $\mathcal S_g$ and are inverse to each other, even though $g$ itself need not be a permutation of $[n]$.
The structural analysis below then shows that the values on the trees of $G(g)$ are forced, while all freedom is concentrated on the cycle core.
The resulting problem on the points lying on directed cycles has a pleasant group-theoretic form when $f$ and $g$ are permutations, and an averaging argument over relabelings extends that viewpoint to arbitrary maps.

\begin{lemma}\label{lem:gfg_change}
The maps $L_g$ and $R_g$ are mutually inverse permutations of $\mathcal S_g$.
\end{lemma}

\begin{proof}
Let $f\in\mathcal S_g$.
For right composition, using $f=g\circ f\circ g$ gives
\[
 g\circ R_g(f)\circ g
 =g\circ f\circ g\circ g
 =f\circ g
 =R_g(f),
\]
so $R_g(f)\in\mathcal S_g$.
Similarly,
\[
 g\circ L_g(f)\circ g
 =g\circ g\circ f\circ g
 =g\circ f
 =L_g(f),
\]
so $L_g(f)\in\mathcal S_g$.
Finally,
\[
 L_g(R_g(f))=g\circ f\circ g=f,
 \qquad
 R_g(L_g(f))=g\circ f\circ g=f.
\]
Thus
\[
 L_g\circ R_g=R_g\circ L_g=\operatorname{id}_{\mathcal S_g},
\]
so $L_g$ and $R_g$ are mutually inverse permutations of $\mathcal S_g$.
\end{proof}

When useful below, we use right composition and write
\[
 h:=R_g(f)=f\circ g.
\]

For permutations, writing $S_n$ for the symmetric group on $[n]$ and $p(n)$ for the number of integer partitions of $n$, we have
\begin{equation}\label{eq:gfg_permutations}
 \#\{(f,g)\in S_n^2:f=g\circ f\circ g\}=n!\,p(n).
\end{equation}
The substitution $u:=f$, $v:=g\circ f$, with inverse $f=u$, $g=v\circ u^{-1}$, transforms the equation into $u^2=v^2$, the problem of Stanley solved by Ehrenborg~\cite{EqualSquares2000}; see also~\cite[Exercises 5.12 and 7.69(h)]{Stanley1999v2}.
Restricting the same count to permutations $g$ with exactly $k$ cycles gives $n!p_k(n)$, where $1\leq k\leq n$, $p_k(n)$ counts partitions into $k$ parts, and fixed points count as cycles of length one.

\subsection{The fixed-map case}

For a fixed map $g$, only its directed cycles affect the number of solutions.

Let $c_m$ denote the number of directed cycles of length $m$ in $G(g)$, and put
\begin{equation}\label{eq:qm_gfg}
 q_m:=\sum_{d\mid m}d\,c_d.
\end{equation}
Thus $q_m$ is the number of vertices in the core $C(g)$ whose cycle length divides $m$.

\begin{theorem}\label{th:fixed_g_gfg}
For a fixed map $g$ on $[n]$, the number of maps $f$ satisfying $f=g\circ f\circ g$ is
\begin{equation}\label{eq:fixed_g_gfg}
 \prod_{m\geq1}q_m^{c_m}
 =
 \prod_{m\geq1}
 \left(\sum_{d\mid m}d\,c_d\right)^{c_m}.
\end{equation}
Factors with $c_m=0$ are omitted from both products.
\end{theorem}

\begin{proof}
By \namedref{Lemma}{lem:gfg_change}, $R_g$ permutes the solution set $\mathcal S_g$.
Thus, for counting purposes, we may replace a solution $f$ by
\[
 h:=R_g(f)=f\circ g,
\]
which is again an element of $\mathcal S_g$ and therefore satisfies
\begin{equation}\label{eq:hghg}
 h=g\circ h\circ g.
\end{equation}

First consider a directed cycle
\[
 C=(x_0,x_1,\dots,x_{m-1}),
 \qquad g(x_i)=x_{i+1},
\]
where subscripts are read modulo $m$.
\namedeqref{Equation}{eq:hghg} gives
\[
 h(x_i)=g(h(x_{i+1})).
\]
Iterating once around $C$ yields
\[
 h(x_0)=g^m(h(x_0)),
\]
so $h(x_0)$ must be a cyclic vertex whose period divides $m$.
There are precisely $q_m$ possibilities.
Once $h(x_0)$ is chosen, the same equations force
$h(x_{m-j})=g^j(h(x_0))$ for $j=1,\dots,m-1$.
Hence each $m$-cycle contributes an independent factor $q_m$.

For the later relabeling argument, we identify these core assignments with maps commuting with the core permutation.
Put $\sigma:=g|_{C(g)}$ and $\eta:=h|_{C(g)}$.
On each $m$-cycle of $\sigma$, use the indexing above and define $\rho(x_i):=x_{-i}$, with subscripts read modulo $m$.
Then $\rho^2=\operatorname{id}$ and
\begin{equation}\label{eq:rho_reverses}
 \rho\circ\sigma\circ\rho=\sigma^{-1}.
\end{equation}
Consequently,
\begin{equation}\label{eq:reverse_to_commute}
 \eta=\sigma\circ\eta\circ\sigma
 \quad\Longleftrightarrow\quad
 \varphi:=\eta\circ\rho\ \text{satisfies}\ 
 \varphi\circ\sigma=\sigma\circ\varphi.
\end{equation}
Indeed, $\eta\circ\sigma^{-1}=\sigma\circ\eta$ and $\rho\circ\sigma=\sigma^{-1}\circ\rho$ give
\[
 \varphi\circ\sigma
 =\eta\circ\sigma^{-1}\circ\rho
 =\sigma\circ\eta\circ\rho
 =\sigma\circ\varphi.
\]
Conversely, if $\varphi$ commutes with $\sigma$, then $\eta:=\varphi\circ\rho$ satisfies
$\sigma\circ\eta\circ\sigma=\varphi\circ\rho=\eta$.
Thus the core contribution also counts the endofunctions commuting with $\sigma$.

There are no further choices on the trees attached to the core.
If $u$ is a noncore vertex and $v=g(u)$ is its parent in $T(g)$, then~\eqref{eq:hghg} gives
\[
 h(u)=g(h(g(u)))=g(h(v)).
\]
Starting at the roots and moving away from the core, this determines $h$ uniquely on every tree vertex.
The resulting map satisfies~\eqref{eq:hghg} at each vertex by construction, so every permitted assignment on the core gives a solution.
Since $R_g$ is a permutation of $\mathcal S_g$, the same count applies to the original solutions $f$.
\end{proof}

\namedref{Figure}{fig:forced-extension} illustrates this extension directly for $f$.
Once $f(1)=1$ is chosen, the cycle equation forces $f(2)=3$ and $f(3)=2$.
The two remaining values follow successively from $f(4)=g(f(3))$ and $f(5)=g(f(4))$.

\begin{figure}[htbp]
\centering
\begin{tikzpicture}[>=stealth,
  point/.style={circle,draw,minimum size=7mm,inner sep=0pt}]
  \node[point] (v1) at (-1.4,0) {$1$};
  \node[point] (v2) at (1.4,0) {$2$};
  \node[point] (v3) at (0,-1.15) {$3$};
  \node[point] (v4) at (0,-2.25) {$4$};
  \node[point] (v5) at (0,-3.35) {$5$};
  \draw[->,thick] (v1) -- (v2);
  \draw[->,thick] (v2) -- (v3);
  \draw[->,thick] (v3) -- (v1);
  \draw[->] (v4) -- (v3);
  \draw[->] (v5) -- (v4);
  \node[above=3mm] at (v1.north) {$f(1)=1$};
  \node[above=3mm] at (v2.north) {$f(2)=3$};
  \node[right=5mm] at (v3.east) {$f(3)=2$};
  \node[right=5mm] at (v4.east) {$f(4)=3$};
  \node[right=5mm] at (v5.east) {$f(5)=1$};
\end{tikzpicture}
\caption{A map $g$ on $[5]$ and one solution of $f=gfg$. The arrows represent $g$, with the cycle drawn thick; the annotations give $f$. The attached-tree values are forced by the values on the cycle.}
\label{fig:forced-extension}
\end{figure}
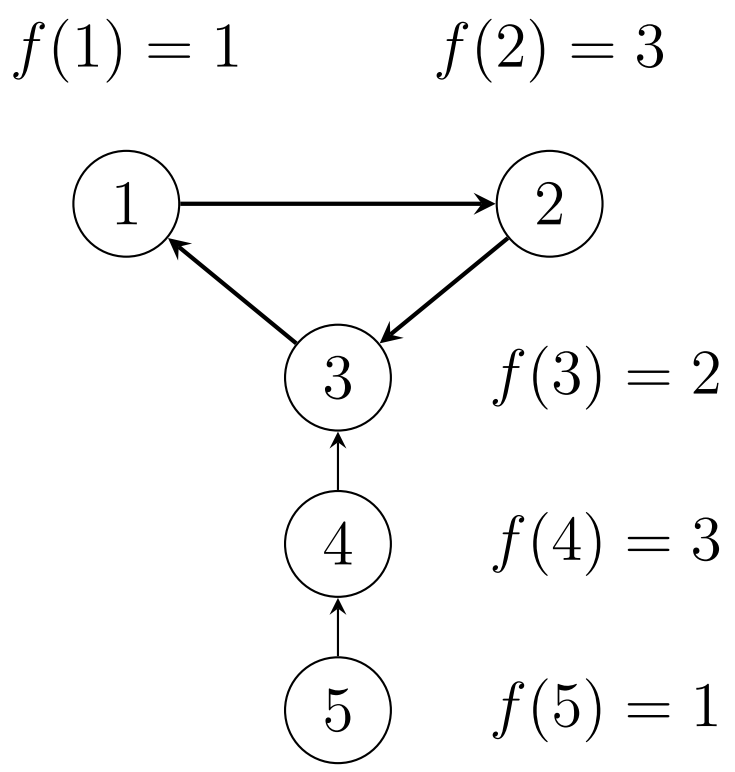

The last proof contains the main structural point of the section: the attached trees of $g$ carry no choices at all.
All freedom lies in choosing, on the permutation core, an endofunction satisfying the corresponding commutation relation.
The product counting maps that commute with a permutation is classical; it appears, for example, in the proof of~\cite[Proposition 2.6]{Holloway2013} and in~\cite[Theorem 2.18(3)]{Daugulis2014}.
The preceding argument explains why the same product applies here even when $g$ has noncyclic vertices.

Every solution takes values in $C(g)$: iterating the equation gives $f=g^r\circ f\circ g^r$, and $g^r([n])=C(g)$ for sufficiently large $r$.
With $\sigma:=g|_{C(g)}$, the equation becomes
\[
 f\circ g=\sigma^{-1}\circ f.
\]
Thus $f$ preserves directed edges from $G(g)$ to $G(\sigma^{-1})$.

\subsection{Averaging over relabelings}

We use \emph{Burnside's lemma}, also called the orbit-counting lemma: if a finite group acts on a finite set, then the number of orbits is the average, over the group elements, of the number of objects they fix~\cite[Lemma 7.24.5]{Stanley1999v2}.

For $k\geq1$, define
\begin{equation}\label{eq:Pk_gfg}
 P_k:=
 \sum_{\substack{c_1,\dots,c_k\geq0\\
                   \sum_{m=1}^k m c_m=k}}
 \frac{k!}{\displaystyle\prod_{m=1}^k m^{c_m}c_m!}
 \prod_{m=1}^k
 \left(\sum_{d\mid m}d\,c_d\right)^{c_m}.
\end{equation}
Factors with $c_m=0$ are omitted from the second product.
The factor preceding the product is the number of permutations on a fixed $k$-element set having permutation cycle type
$1^{c_1}2^{c_2}\cdots k^{c_k}$, meaning that there are $c_m$ cycles of length $m$.
Thus $P_k$ is the total core contribution after summing \namedref{Theorem}{th:fixed_g_gfg} over all core permutations.

By the notation introduced before \namedref{Theorem}{th:number_gfg}, $a_k$ is the number of endofunctions on a $k$-element set up to relabeling, equivalently the number of unlabeled functional digraphs on $k$ vertices.

\begin{proposition}\label{prop:Pk_unlabeled}
For every $k\geq1$,
\begin{equation}\label{eq:Pk_ak}
 P_k=k!\,a_k.
\end{equation}
\end{proposition}

\begin{proof}
The symmetric group $S_k$ acts on maps $\varphi:[k]\to[k]$ by relabeling. For $\sigma\in S_k$, the action is
\[
 \varphi\longmapsto \sigma\circ\varphi\circ\sigma^{-1}.
\]
The orbits, i.e., the classes obtained by relabeling, are precisely the unlabeled endofunctions.
An endofunction is fixed by $\sigma$ exactly when it commutes with $\sigma$.
If $\sigma$ has $c_m$ cycles of length $m$, the cycle calculation in the proof of \namedref{Theorem}{th:fixed_g_gfg} shows that the number of such commuting endofunctions is
\[
 \prod_{m\geq1}
 \left(\sum_{d\mid m}d\,c_d\right)^{c_m}.
\]
Burnside's lemma therefore gives
\[
 a_k=\frac1{k!}\sum_{\sigma\in S_k}
 \#\{\varphi:\varphi\circ\sigma=\sigma\circ\varphi\}.
\]
Grouping this sum by permutation cycle type gives exactly~\eqref{eq:Pk_gfg}, proving~\eqref{eq:Pk_ak}.
\end{proof}

The classical numbers of unlabeled functional digraphs begin as follows~\cite{Harary1959}:
\[
 a_0,a_1,a_2,a_3,a_4,a_5,a_6
 =1,1,3,7,19,47,130.
\]
Thus
\[
 P_1,P_2,P_3,P_4,P_5,P_6
 =1,6,42,456,5640,93600.
\]

In the permutation count~\eqref{eq:gfg_permutations}, the factor $p(n)$ counts permutations up to relabeling.
Their ordinary generating function is the classical partition series
\begin{equation}\label{eq:partition_ogf}
 \sum_{k\geq0}p(k)u^k
 =\prod_{j\geq1}\frac1{1-u^j}.
\end{equation}
For arbitrary endofunctions, each cycle vertex carries an unlabeled rooted tree.
The decomposition recorded in~\eqref{eq:A_R_product} therefore replaces the bare cycle contribution $u^j$ in~\eqref{eq:partition_ogf} by the rooted-tree contribution $\mathcal R(u^j)$; the rooted-tree series itself satisfies~\eqref{eq:Polya_rooted_tree}.
This is the unlabeled-tree layer of the problem.

\subsection{Attaching labeled trees: the composition \texorpdfstring{$\mathcal A\circ\Tree$}{A o T}}

We now restore the noncyclic vertices of $g$.  These attached trees are labeled, so the labeled tree function
\[
 \Tree(x)=x e^{\Tree(x)}=-W(-x)
\]
from \namedref{Section}{sec:graphs} enters for a third time.

We are now ready to prove the third main theorem.

\begin{proof}[Proof of \namedref{Theorem}{th:number_gfg}]
Fix the number $k$ of cyclic vertices of $g$.
Choose them in $\binom{n}{k}$ ways.
After summing over the possible core permutations, \namedref{Proposition}{prop:Pk_unlabeled} gives the core factor $P_k=k!a_k$.
It remains to attach the $n-k$ noncyclic vertices to the prescribed $k$ core vertices.
The rooted-forest form of Cayley's formula gives
\[
 k n^{n-k-1}
\]
such forests, with all arcs directed toward their specified roots; see, for example,~\cite[Note II.18, p.~128]{FlajoletSedgewick2009}.
This proves~\eqref{eq:Dn_gfg}.

For the generating-function form, Lagrange inversion applied to
$\Tree(x)=x e^{\Tree(x)}$ gives, for $1\leq k\leq n$,
\begin{equation}\label{eq:Tk_lagrange}
 \coef{x^n}\Tree(x)^k
 =\frac{k}{n}\coef{u^{n-k}}e^{nu}
 =\frac{k n^{n-k-1}}{(n-k)!};
\end{equation}
see~\cite{SuryaWarnke2023}.
Hence
\begin{align*}
 n!\coef{x^n}\mathcal A(\Tree(x))
 &=\sum_{k=1}^n a_k\,n!\coef{x^n}\Tree(x)^k\\
 &=\sum_{k=1}^n
   \binom nk\,k n^{n-k-1}k!a_k
 =D_n.
\end{align*}
The constant terms agree as well, since $D_0=a_0=1$.
By~\eqref{eq:tree_W}, the composition in~\eqref{eq:D_composition} can equivalently be written as
\[
 \mathcal D(x)=\mathcal A(-W(-x)).
\]
Finally, substituting $u=\Tree(x)$ into~\eqref{eq:A_R_product} gives the useful two-tree form
\begin{equation}\label{eq:D_two_tree_functions}
 \mathcal D(x)=\prod_{j\geq1}\frac{1}{1-\mathcal R(\Tree(x)^j)}.
\end{equation}
\end{proof}

Thus \namedref{Section}{sec:gfg} has two closely related tree layers.
The Burnside average on the core produces the ordinary generating function $\mathcal A(u)$ of \emph{unlabeled} functional digraphs, whose cycle-and-tree decomposition is governed by the unlabeled rooted-tree series $\mathcal R(u)$.
The vertices outside the core form \emph{labeled} rooted trees, whose exponential generating function is $\Tree(x)=-W(-x)$.
The compact composition $\mathcal D=\mathcal A\circ\Tree$ records both layers at once.

Together with $D_0=1$, the first values are
\[
 1,\ 1,\ 8,\ 87,\ 1312,\ 24685,\ 564336,\ 15168811,\dots.
\]
They form OEIS sequence A239753~\cite{OEIS}, corresponding to
$f(x)=g(f(g(x)))$.
The agreement provides a numerical check on \namedref{Theorem}{th:number_gfg}.

The permutation subcase and the full mapping problem now sit in the same framework.
For permutations, the core is the whole set and the count is $n!p(n)$.
For arbitrary maps, Burnside's lemma averages the commuting endofunctions on the core, and the noncyclic vertices contribute the labeled-tree substitution $\Tree(x)=-W(-x)$.
We next apply this fixed-map counting viewpoint to commutation with multiplication on residue classes.
The arithmetic of the multiplier will determine both the directed cycles and the choices on the attached trees.

\section{The identity \texorpdfstring{$f(ax)=af(x)$}{f(ax) = a f(x)} modulo \texorpdfstring{$n$}{n}}\label{sec:px}

The preceding sections first fixed one map and counted the compatible choices of the other.
We now carry out that count for a prescribed arithmetic family.
Let $a$ be a positive integer and define
\[
 g_a(x):=ax\pmod n
\]
on $\mathbb Z_n:=\mathbb Z/n\mathbb Z$.
Then
\begin{equation}\label{eq:ax}
 f(ax)=a f(x)
\end{equation}
is simply the commutation relation $f\circ g_a=g_a\circ f$, with $g_a$ fixed.
Functional digraphs of multiplication maps on finite cyclic groups have been studied in considerable generality; see, for example, Qureshi and Panario~\cite{QureshiPanario2015}.
Higgins~\cite[Theorem~5]{Higgins1988} gives a functional-digraph construction for maps $q$ satisfying $q\circ s=t\circ q$ for prescribed maps $s,t$; commutation is the case $s=t$.
Ara\'ujo and Konieczny~\cite[Corollary~5.1]{AraujoKonieczny2013} describe all maps commuting with a fixed finite map in terms of its directed cycles and attached branches; see also Daugulis~\cite{Daugulis2014}.

Our aim here is more elementary and explicit. We first treat a prime multiplier. In that case every tree attached to a point on a directed cycle has the same regular branching structure, and the count of commuting maps collapses to a simple product. We then identify a composite setting in which exactly the same regularity survives.

\subsection{Prime multipliers}

Let $p$ be a fixed prime, and let $\nu_p(n)$ denote the exponent of $p$ in the prime factorization of $n$.
For a residue class $k\in\mathbb Z_n$, divisibility statements such as $p^j\mid k$ mean that the class has an integer representative divisible by $p^j$; for $j\leq\nu_p(n)$ this is independent of the chosen representative.
For $q>1$ and $\gcd(b,q)=1$, $\ord_q(b)$ denotes the \emph{multiplicative order} of $b$ modulo $q$, i.e., the least positive $r$ for which $b^r\equiv1\pmod q$; we use the convention $\ord_1(b):=1$.
We also use the arithmetic M\"obius function $\mu$ introduced before \namedref{Theorem}{th:number_px}.
The divisor-sum inversion rule called \emph{M\"obius inversion} says that if $F(\ell)=\sum_{d\mid\ell}G(d)$, then
\[
 G(\ell)=\sum_{d\mid\ell}\mu(\ell/d)F(d).
\]

Write
\[
 n=p^tm,\qquad t:=\nu_p(n),\qquad p\nmid m,
\]
and let $g(x):=px$ on $\mathbb Z_n$.
The \emph{depth} of a vertex in $T(g)$ is its distance from the root of its tree.
We begin by describing $G(g)$.

\begin{theorem}\label{th:Gpx}
The functional digraph $G(g)$ has the following properties.
\begin{enumerate}
\item If $t\geq1$, a vertex $k$ divisible by $p$ has exactly $p$ incoming arcs, namely
\[
 \left(k/p+s\,n/p,\,k\right),\qquad s=0,1,\dots,p-1,
\]
whereas a vertex not divisible by $p$ has no incoming arcs.

\item If $t=0$, multiplication by $p$ is a permutation of $\mathbb Z_n$, and every vertex $k$ has the unique incoming arc $(p^{-1}k,k)$.

\item A vertex $k$ lies on a directed cycle if and only if $p^t\mid k$.

\item If $k$ is cyclic, then the length of its cycle is
\[
 \ord_{m/\gcd(m,k)}(p).
\]

\item The core $C(g)$ contains exactly $m$ vertices, and hence $T(g)$ has exactly $m$ rooted trees.

\item The number $N_\ell$ of vertices lying on cycles of length exactly $\ell$ is
\begin{equation}\label{eq:Nell}
 N_\ell=\sum_{d\mid\ell}\mu(\ell/d)\gcd(m,p^d-1).
\end{equation}

\item If $t\geq1$, every tree in $T(g)$ has the same shape: its root has $p-1$ children, and from each of these children grows a complete $p$-ary tree of height $t-1$, meaning that every nonleaf vertex has exactly $p$ children and all leaves are at the same depth.
\end{enumerate}
\end{theorem}

\begin{proof}
For (1), suppose first that $t\geq1$.
If $p\mid k$, the congruence
$pu\equiv k\pmod n$ has the $p$ solutions
\[
 u\equiv k/p+s\,n/p\pmod n,
 \qquad s=0,1,\dots,p-1.
\]
If $p\nmid k$, it has no solution already modulo $p$.
This proves (1).
Statement (2) follows because $p$ is invertible modulo $n$ when $t=0$.

For (3), suppose $k$ lies on a cycle.
Then $p^\ell k\equiv k\pmod n$ for some $\ell>0$, so
\[
 (p^\ell-1)k\equiv0\pmod{p^tm}.
\]
Because $p\nmid p^\ell-1$, this implies $p^t\mid k$.
Conversely, if $k=p^ta$, then multiplication by $p$ is invertible modulo $m$, so some positive power of $p$ fixes $a$ modulo $m$.
The corresponding power of $g$ fixes $k$ modulo $n$, proving that $k$ is cyclic.

For (4), write the cyclic vertex as $k=p^ta$. Then $p^\ell k\equiv k\pmod n$ is equivalent to
\[
 m\mid a(p^\ell-1),
\]
or, equivalently,
\[
 p^\ell\equiv1\pmod{m/\gcd(m,a)}.
\]
Since $p^t$ is coprime to $m$, we have $\gcd(m,a)=\gcd(m,k)$.
The least positive such $\ell$ is the stated multiplicative order.

The cyclic vertices are exactly the classes $p^ta$ with $0\leq a<m$, so there are $m$ of them; this proves (5).

For (6), the cyclic vertices whose periods divide $\ell$ are precisely the solutions of $g^\ell(k)=k$ in the core.
Writing $k=p^t x$, this condition becomes
\[
 (p^\ell-1)x\equiv0\pmod m,
\]
which has $\gcd(m,p^\ell-1)$ solutions.
Thus
\[
 \gcd(m,p^\ell-1)=\sum_{d\mid\ell}N_d,
\]
and the M\"obius inversion rule above gives~\eqref{eq:Nell}.

Finally, let $t\geq1$ and consider a root $p^ta$.
Its $p$ preimages have the form
\[
 p^{t-1}(a+ms),\qquad s=0,1,\dots,p-1.
\]
Exactly one choice of $s$ satisfies
\[
 s\equiv-a m^{-1}\pmod p;
\]
for that choice the preimage is divisible by $p^t$ and is the predecessor of the root on its directed cycle.
The remaining $p-1$ preimages are the children of the root in $T(g)$.

More generally, a vertex at depth $j$, where $1\leq j<t$, is divisible by $p^{t-j}$ but not by $p^{t-j+1}$.
Each of its $p$ preimages is divisible by $p^{t-j-1}$ but not by $p^{t-j}$, so all $p$ preimages lie one level farther from the root.
Vertices at depth $t$ are not divisible by $p$ and have no preimages.
This proves (7).
\end{proof}

We next separate the choices on the core from those in the attached trees. In \namedref{Figure}{fig:forced-extension}, the image of a child was forced to be $g(f(v))$, where $v$ is its parent. For commutation, the child may instead map to any preimage of $f(v)$ under $g$. The regular branching makes the number of choices uniform.

\Needspace{7\baselineskip}
\begin{theorem}\label{th:commuting_extension}\leavevmode
\begin{enumerate}
\item If $f\circ g=g\circ f$, then a cycle of length $L$ in $G(g)$ is mapped by $f$ onto a cycle of length $\ell$ dividing $L$.
For a fixed target cycle of length $\ell\mid L$, there are exactly $\ell$ possible restrictions of $f$ to the domain cycle.

\item Every map $h:C(g)\to C(g)$ commuting with $g|_{C(g)}$ has exactly
\[
 p^{n-m}
\]
extensions $f:\mathbb Z_n\to\mathbb Z_n$ satisfying $f\circ g=g\circ f$.
\end{enumerate}
\end{theorem}

\begin{proof}
Let
\[
 C=(v,g(v),g^2(v),\dots,g^{L-1}(v))
\]
be a cycle with $g^L(v)=v$.
Commutation gives
\[
 f(g^i(v))=g^i(f(v))
\]
for every $i\geq0$.
In particular, $g^L(f(v))=f(v)$, so $f(v)$ lies on a cycle whose length $\ell$ divides $L$.
Once a target cycle of length $\ell$ is fixed, $f(v)$ may be chosen as any of its $\ell$ vertices, and this choice uniquely determines $f$ on all of $C$.
This proves (1).

Now suppose $f$ is fixed on the core.
If $v$ has depth $j$, then $g^j(v)$ lies on the core, and therefore
\[
 g^j(f(v))=f(g^j(v))
\]
also lies on the core.
Thus $f(v)$ has depth at most $j$; in particular, $f$ cannot increase depth.

Suppose $u$ is a child of $v$ in $T(g)$ and $f(v)$ has already been chosen.
Since $u$ is not a root, the depth of $v$ is at most $t-1$.
Hence the depth of $f(v)$ is also at most $t-1$, and by \namedref{Theorem}{th:Gpx} the vertex $f(v)$ has exactly $p$ incoming arcs in $G(g)$.
The commutation relation requires
\[
 g(f(u))=f(g(u))=f(v),
\]
so $f(u)$ may be chosen independently as any of these $p$ preimages.
Every such preimage has depth at most one more than that of $f(v)$, hence at most the depth of $u$.
Thus the same argument applies successively to all descendants, and no choice prevents a later extension.
There are $n-m$ nonroot vertices in $T(g)$, giving $p^{n-m}$ extensions in total.
\end{proof}

Recall the numbers $N_d$ from~\eqref{eq:Nd}.  By \namedref{Theorem}{th:Gpx}, $N_d$ is the number of core vertices of exact period $d$, and $N_d/d$ is the number of directed cycles of length $d$ in $G(g)$.

We can now prove \namedref{Theorem}{th:number_px}.

\begin{proof}[Proof of \namedref{Theorem}{th:number_px}]
Assume first that $m>1$.
Every directed cycle in the core of $G(g)$ has length dividing $r=\ord_m(p)$.
Fix a directed domain cycle of length $d$ in $G(g)$.
By \namedref{Theorem}{th:commuting_extension}, it can be mapped to any target cycle in $G(g)$ whose length $\ell$ divides $d$.
There are $N_\ell/\ell$ target cycles of length $\ell$ in $G(g)$, and each target cycle allows $\ell$ choices for the image of one selected vertex of the domain cycle.
Hence the number of possible images of a fixed $d$-cycle is
\[
 \sum_{\ell\mid d}\frac{N_\ell}{\ell}\,\ell
 =\sum_{\ell\mid d}N_\ell
 =\gcd(m,p^d-1),
\]
where the last equality is the relation inverted in~\eqref{eq:Nd}.

There are $N_d/d$ domain cycles of length $d$, and their image choices are independent.
Therefore the number of possible restrictions of $f$ to the core is
\[
 \prod_{d\mid r}\gcd(m,p^d-1)^{N_d/d}.
\]
By \namedref{Theorem}{th:commuting_extension}, each such restriction has $p^{n-m}$ extensions to the whole functional digraph, proving~\eqref{eq:Bpn}.

If $m=1$, the core consists only of the fixed point $0$, which must map to itself.
\namedref{Theorem}{th:commuting_extension} then gives $p^{n-1}$ extensions.
\end{proof}

\begin{corollary}\label{cor:p2}
If $n=2^tm$ with $m$ odd and $m>1$, then the number of maps $f$ on $\mathbb Z_n$ satisfying $f(2x)=2f(x)$ is
\[
 2^{n-m}
 \prod_{d\mid\ord_m(2)}
 \gcd(m,2^d-1)^{\frac{1}{d}
 \sum_{q\mid d}\mu(d/q)\gcd(m,2^q-1)}.
\]
If $m=1$, the number is $2^{n-1}$.
\end{corollary}

For $n=1,2,3,\dots$, these numbers begin
\[
 1,\ 2,\ 3,\ 8,\ 5,\ 24,\ 49,\ 128,\ 27,\ 160,\dots,
\]
and form OEIS sequence A117987~\cite{OEIS}.

\subsection{A composite extension}

The prime hypothesis in the preceding argument is stronger than necessary. What matters for the extension count is the uniform shape of the trees attached to the core. The same shape occurs for a natural family of composite multipliers, which we now define.

Let $a$ be \emph{squarefree}, i.e., divisible by no square greater than $1$, and put $d:=\gcd(a,n)$. We say that the pair $(a,n)$ is \emph{balanced} if either $d=1$, or there is an integer $t\geq1$ such that
\begin{equation}\label{eq:balanced}
 n=d^t m,\qquad \gcd(a,m)=1.
\end{equation}
Equivalently, when $d>1$, every prime dividing both $a$ and $n$ occurs in $n$ with the same exponent $t$.

The case $a=p$ prime is automatically balanced: then $d=p$ whenever $p\mid n$, and~\eqref{eq:balanced} is simply the usual factorization $n=p^t m$ with $p\nmid m$.

For the next structural result, assume $d>1$ and use the factorization
\[
 n=d^t m,\qquad t\geq1,\qquad \gcd(a,m)=1
\]
from~\eqref{eq:balanced}; let $g_a(x):=ax$ on $\mathbb Z_n$.

\Needspace{6\baselineskip}
\begin{theorem}\label{th:balanced_graph}
The functional digraph $G(g_a)$ has the following properties.
\begin{enumerate}
\item A vertex $k$ has $d$ incoming arcs if $d\mid k$, and no incoming arcs otherwise.

\item A vertex $k$ is cyclic if and only if $d^t\mid k$. Consequently the core has exactly $m$ vertices.

\item If $k$ is cyclic, its cycle length is
\[
 \ord_{m/\gcd(m,k)}(a).
\]

\item The number $M_\ell$ of core vertices of exact period $\ell$ is
\begin{equation}\label{eq:Mell}
 M_\ell=\sum_{q\mid\ell}\mu(\ell/q)\gcd(m,a^q-1).
\end{equation}

\item Every tree in $T(g_a)$ has the same shape: its root has $d-1$ children, and from each of them grows a complete $d$-ary tree of height $t-1$ in the same sense as the complete $p$-ary trees above.
\end{enumerate}
\end{theorem}

\begin{proof}
Write $a=du$. Since $a$ is squarefree and $d$ contains precisely the prime factors of $a$ that divide $n$, we have $\gcd(u,n)=1$.

The congruence $ax\equiv k\pmod n$ has solutions if and only if $d\mid k$; when it does, it has exactly $d$ solutions. This proves (1).

If $k$ is cyclic, then $a^\ell k\equiv k\pmod n$ for some $\ell>0$. For every prime $q\mid d$ we have $q\mid a$ and hence $q\nmid a^\ell-1$. Since $q^t\mid n$, it follows that $q^t\mid k$. As $d$ is squarefree, this holds simultaneously for all $q\mid d$ exactly when $d^t\mid k$. Conversely, if $k=d^t b$, then multiplication by $a$ is invertible modulo $m$, so some positive power of $a$ fixes $b$ modulo $m$ and therefore fixes $k$ modulo $n$. This proves (2).

For a cyclic vertex $k=d^tb$, the congruence $a^\ell k\equiv k\pmod n$ reduces to
\[
 (a^\ell-1)b\equiv0\pmod m.
\]
Since $d$ is coprime to $m$, $\gcd(m,b)=\gcd(m,k)$, and (3) follows.

The number of core vertices fixed by $g_a^\ell$ is the number of solutions of
\[
 (a^\ell-1)b\equiv0\pmod m,
\]
namely $\gcd(m,a^\ell-1)$. The M\"obius inversion rule stated above gives (4).

It remains to describe the attached trees. A root is divisible by $d^t$ and therefore has $d$ preimages. Multiplication by $a$ restricts to a permutation of the core, so exactly one of these preimages belongs to the core. That preimage is the predecessor on the cycle, leaving $d-1$ children in the tree.

More generally, let $v$ be a noncyclic vertex at depth $j$, where $1\leq j<t$.
Since $a^jv$ lies on the core, $d^t\mid a^jv$.
Squarefreeness of $a$ therefore gives $d^{t-j}\mid v$, so $v$ has exactly $d$ preimages.
Every such preimage has depth $j+1$. Also, $d^t\mid a^t$, so every vertex reaches the core after at most $t$ iterations. Hence a vertex at depth $t$ cannot have a preimage and, by (1), is not divisible by $d$. This gives the stated $d$-ary structure.
\end{proof}

For example, take $a=6$ and $n=36=6^2$. The core consists only of $0$. After its loop is deleted, the root has $5$ children, and each of those has $6$ children, accounting for all $1+5+30=36$ vertices. This is the same regular geometry as in the prime case, with the branching number $6$ in place of $p$.

The preceding proof shows why the prime argument extends without any additional counting machinery: one simply replaces the branching number $p$ by $d$.

For the enumeration, let $a$ be squarefree and let $(a,n)$ be balanced.
Put $d:=\gcd(a,n)$, and write $n=d^t m$ with $\gcd(a,m)=1$ when $d>1$; if $d=1$, put $m:=n$.
For $e\geq1$, define
\begin{equation}\label{eq:Md}
 M_e:=\sum_{q\mid e}\mu(e/q)\gcd(m,a^q-1).
\end{equation}
If $m>1$, let $r:=\ord_m(a)$.

\begin{theorem}\label{th:number_balanced}
If $m>1$, the number $B_a(n)$ of maps $f$ on $\mathbb Z_n$ satisfying $f(ax)=a f(x)$ is
\begin{equation}\label{eq:Ban}
 B_a(n)=d^{n-m}
 \prod_{e\mid r}\gcd(m,a^e-1)^{M_e/e}.
\end{equation}
When $m=1$, the number is $d^{n-1}$. For $d=1$, the factor $d^{n-m}$ is understood as $1$.
\end{theorem}

\begin{proof}
If $d=1$, multiplication by $a$ is a permutation, so there are no nonroot vertices and only the cycle calculation is needed. Assume therefore that $d>1$.
By \namedref{Theorem}{th:balanced_graph}, the core has $m$ vertices, and its cycle counts are given by~\eqref{eq:Md}. Exactly as in the proof of \namedref{Theorem}{th:number_px}, a directed domain cycle of length $e$ in $G(g_a)$ has
\[
 \sum_{q\mid e}M_q=\gcd(m,a^e-1)
\]
possible images on the core. Multiplying over the $M_e/e$ domain cycles of each length gives the product in~\eqref{eq:Ban}.

It remains to extend a commuting map from the core to the attached trees. The depth argument from \namedref{Theorem}{th:commuting_extension} applies verbatim. By \namedref{Theorem}{th:balanced_graph}, every possible image of the parent of a nonroot vertex has exactly $d$ preimages. Thus each of the $n-m$ nonroot vertices contributes an independent factor $d$, giving $d^{n-m}$ extensions. If $m=1$, the core is the single fixed point $0$, and the same argument gives $d^{n-1}$.
\end{proof}

The prime theorem is recovered from \namedref{Theorem}{th:number_balanced} by taking $a=p$, in which case $d=p$ whenever $p\mid n$. Thus primality is not itself responsible for the simple formula; the essential feature is the regularity of the trees attached outside the directed cycles.

\medskip\noindent\textit{Why arbitrary multipliers are less uniform.}\ 

Squarefreeness alone does not guarantee the regular tree structure. For instance, take $a=6$ and $n=18$. Then $d=\gcd(6,18)=6$, but $18$ cannot be written as $6^t m$ with $\gcd(6,m)=1$. The preimages of $0$ are
\[
 0,3,6,9,12,15.
\]
After the loop at $0$ is removed, the root has five children. However, among these children only $6$ and $12$ are divisible by $6$ and therefore have preimages; $3,9,$ and $15$ are already leaves. The uniform $d$-ary tree from \namedref{Theorem}{th:balanced_graph} has disappeared.

For a completely general multiplier $a$, the successive layers are controlled by the kernels of the iterates of $g_a$, where
$\ker(g_a^j):=\{x\in\mathbb Z_n:g_a^j(x)=0\}$; in particular,
\[
 |\ker(g_a^j)|=\gcd(n,a^j).
\]
This describes the possible \emph{preperiods}, i.e., the numbers of steps needed before an orbit first reaches a directed cycle, but the number of choices for the image of a nonperiodic vertex now depends on where its parent was mapped.
The enumeration therefore becomes a problem of graph homomorphisms of the full functional digraph, i.e., maps of vertices that preserve directed arcs.
Such a map may send an attached-tree edge to a cycle edge, and its enumeration requires the shapes of the attached rooted trees rather than a single independent branching factor.
General descriptions of multiplication-map trees appear in~\cite{QureshiPanario2015}, and commuting maps are treated as digraph homomorphisms in~\cite{Daugulis2014}.

The prime and balanced squarefree cases make the contribution of the trees explicit: after choosing a compatible map on the cycles, every noncyclic vertex supplies the same number of choices.

\section{Concluding remarks}\label{sec:conclusion}

The four problems share a fixed-map counting step, but the information needed in that step differs.
For $f\circ g=f\circ g\circ g$, sources and components suffice.
For $f\circ f=f\circ g\circ f$, the useful data are the fibers of $f$ and the positions of its image within those fibers.
For $f=g\circ f\circ g$, the restriction to the cycle core determines all values on the attached trees, and averaging over relabelings introduces unlabeled endofunctions.
The arithmetic application retains the fixed-map viewpoint: the regular trees of a multiplication map supply an explicit extension factor.

\section*{Acknowledgments}

ChatGPT assisted with developing and reviewing mathematical arguments, checking calculations, and editing the manuscript. The author accepts full responsibility for its contents.

\bibliographystyle{plainurl}
\bibliography{refs_rev131}

\begin{thebibliography}{10}

\bibitem{AlekseyevType23}
Max~A. Alekseyev.
\newblock Enumerating solutions to type {(2,3)} composition identities of
  finite maps.
\newblock Unpublished follow-up manuscript, 2026.

\bibitem{AraujoKonieczny2013}
Jo{\~a}o Ara{\'u}jo and Janusz Konieczny.
\newblock Centralizers in the full transformation semigroup.
\newblock {\em Semigroup Forum}, 86(1):1--31, 2013.
\newblock \href {https://doi.org/10.1007/s00233-012-9424-0}
  {\path{doi:10.1007/s00233-012-9424-0}}.

\bibitem{Corless1996}
Robert~M. Corless, Gaston~H. Gonnet, D.~E.~G. Hare, David~J. Jeffrey, and
  Donald~E. Knuth.
\newblock On the {L}ambert {$W$} function.
\newblock {\em Advances in Computational Mathematics}, 5:329--359, 1996.
\newblock \href {https://doi.org/10.1007/BF02124750}
  {\path{doi:10.1007/BF02124750}}.

\bibitem{Daugulis2014}
Peteris Daugulis.
\newblock Graph structure of commuting functions, 2014.
\newblock \href {https://arxiv.org/abs/1407.6035} {\path{arXiv:1407.6035}}.

\bibitem{EqualSquares2000}
Richard Ehrenborg.
\newblock Pairs with equal squares: 10654.
\newblock {\em The American Mathematical Monthly}, 107(4):368, 2000.
\newblock Solution to a problem proposed by Richard P. Stanley.
\newblock \href {https://doi.org/10.2307/2589191} {\path{doi:10.2307/2589191}}.

\bibitem{EinsteinNotebooks}
David Einstein.
\newblock {oeis-function-identities}.
\newblock Mathematical notebooks: EasyCounts; Pseudoinverses; Cyclic Count;
  Generalized tree polynomials. Accessed November 3, 2026.
\newblock URL: \url{https://github.com/deinst/oeis-function-identities}.

\bibitem{FlajoletSedgewick2009}
Philippe Flajolet and Robert Sedgewick.
\newblock {\em Analytic Combinatorics}.
\newblock Cambridge University Press, Cambridge, 2009.
\newblock \href {https://doi.org/10.1017/CBO9780511801655}
  {\path{doi:10.1017/CBO9780511801655}}.

\bibitem{Harary1959}
Frank Harary.
\newblock The number of functional digraphs.
\newblock {\em Mathematische Annalen}, 138:203--210, 1959.
\newblock \href {https://doi.org/10.1007/BF01342903}
  {\path{doi:10.1007/BF01342903}}.

\bibitem{Higgins1988}
Peter~M. Higgins.
\newblock Digraphs and the semigroup of all functions on a finite set.
\newblock {\em Glasgow Mathematical Journal}, 30(1):41--57, 1988.
\newblock \href {https://doi.org/10.1017/S0017089500007011}
  {\path{doi:10.1017/S0017089500007011}}.

\bibitem{Holloway2013}
Michael Holloway and Mark Shattuck.
\newblock Commuting pairs of functions on a finite set.
\newblock {\em Pure Mathematics and Applications}, 24(1):9--18, 2013.

\bibitem{OEIS}
{OEIS Foundation Inc.}
\newblock The on-line encyclopedia of integer sequences.
\newblock Numerical data for the sequences cited in this article; accessed
  November 3, 2026.
\newblock URL: \url{https://oeis.org/}.

\bibitem{QureshiPanario2015}
Claudio Qureshi and Daniel Panario.
\newblock {R}\'edei actions on finite fields and multiplication map in cyclic
  group.
\newblock {\em SIAM Journal on Discrete Mathematics}, 29(3):1486--1503, 2015.
\newblock \href {https://doi.org/10.1137/140993338}
  {\path{doi:10.1137/140993338}}.

\bibitem{Stanley1999v2}
Richard~P. Stanley.
\newblock {\em Enumerative Combinatorics}, volume~2.
\newblock Cambridge University Press, Cambridge, 1999.
\newblock \href {https://doi.org/10.1017/CBO9780511609589}
  {\path{doi:10.1017/CBO9780511609589}}.

\bibitem{SuryaWarnke2023}
Erlang Surya and Lutz Warnke.
\newblock Lagrange inversion formula by induction.
\newblock {\em The American Mathematical Monthly}, 130(10):944--948, 2023.
\newblock \href {https://doi.org/10.1080/00029890.2023.2251344}
  {\path{doi:10.1080/00029890.2023.2251344}}.

\end{thebibliography}

\clearpage
\appendix
\refstepcounter{section}
\section*{Appendix \thesection: Eliminating the constant term}\label{sec:Cn_univariate}

Fix $n$ and use the series $Q_n(y,z)$ from~\eqref{eq:Qn}. To replace the constant term in~\eqref{eq:Cn_constant_Q} by a one-variable coefficient, we solve $z=Q_n(y,z)$ as a formal series in $y$. The auxiliary series below is the formal root with zero constant term.
Define $S_n(y)\in y\mathbb Q[[y]]$ uniquely by
\begin{equation}\label{eq:Sn}
 S_n(y)=Q_n\bigl(y,S_n(y)\bigr).
\end{equation}
Since $Q_n(y,z)$ is divisible by $y$, the coefficient of $y^r$ on the right depends only on coefficients of $S_n(y)$ of degree less than $r$.
Thus the coefficients of $S_n(y)$ are determined successively, proving existence and uniqueness.
A constant-term form of Lagrange inversion gives
\begin{equation}\label{eq:Cn_Sn}
 C_n=n!\,\coef{y^n}
 \left(1-Q_{n,z}\bigl(y,S_n(y)\bigr)\right)^{-1},
\end{equation}
where $Q_{n,z}$ denotes the partial derivative with respect to $z$.
For completeness, we justify this identity in the formal series ring already specified.
Write $Q:=Q_n$ and $S:=S_n(y)$, and define
\[
 E(y,z):=\frac{Q(y,z)-Q(y,S)}{z-S}.
\]
The division is exact coefficientwise in $y$, since each coefficient of $Q$ is a polynomial in its second argument.
Hence $E\in y\mathbb Q[z][[y]]$, and $E(y,S)=Q_{n,z}(y,S)$.
Using $Q(y,S)=S$, we obtain
\[
 z-Q(y,z)=(z-S)(1-E(y,z)).
\]
For any $H\in\mathbb Q[z][[y]]$, the geometric expansion gives
\[
 \coef{z^0}\frac{z}{z-S}H(y,z)
 =\coef{z^0}\sum_{j\geq0}S^jz^{-j}H(y,z)
 =H(y,S).
\]
All sums are well-defined coefficientwise because $S\in y\mathbb Q[[y]]$.
Taking $H=(1-E)^{-1}$ proves~\eqref{eq:Cn_Sn}.
Thus~\eqref{eq:Cn_Sn} gives a one-variable coefficient formula for $C_n$; the series being used depends on $n$.

The tree form~\eqref{eq:Qn_tree} turns this immediately into the following explicit version of the one-variable formula. Set $U_n(y):=\Tree(yS_n(y))$.

\begin{corollary}\label{cor:Cn_lagrange}
For the series $S_n(y)$ defined by~\eqref{eq:Sn},
\begin{equation}\label{eq:Cn_lagrange}
 C_n=n!\,\coef{y^n}
 \left(
  1-
  \frac{
   U_n(y)\bigl(1+S_n(y)\bigr)
   \left(S_n(y)+nU_n(y)
   \left(1-U_n(y)\right)\right)
  }{
   S_n(y)^2\left(1-U_n(y)\right)^2
  }
 \right)^{-1}.
\end{equation}
\end{corollary}

\Needspace{10\baselineskip}
\begin{proof}
Write $T:=\Tree(yz)$.  Differentiating~\eqref{eq:Qn_tree} with respect to $z$, and using
\[
 T_z=\frac{T}{z(1-T)},
\]
gives
\begin{equation}\label{eq:Qnz_tree}
 Q_{n,z}(y,z)
 =\bigl(1+Q_n(y,z)\bigr)
  \frac{T\bigl(z+nT(1-T)\bigr)}
       {z^2(1-T)^2}.
\end{equation}
Now set $z=S_n(y)$ and use $Q_n(y,S_n(y))=S_n(y)$.  This gives
\[
 Q_{n,z}\bigl(y,S_n(y)\bigr)
 =\frac{
   U_n(y)\bigl(1+S_n(y)\bigr)
   \left(S_n(y)+nU_n(y)
   \left(1-U_n(y)\right)\right)
  }{
   S_n(y)^2\left(1-U_n(y)\right)^2
  }.
\]
Although displayed with $S_n(y)^2$ in the denominator, this is a formal power series: indeed,
\[
 \frac{\Tree(yS_n(y))}{S_n(y)}
 =y\,\frac{\Tree(yS_n(y))}{yS_n(y)}\in y\mathbb Q[[y]].
\]
Substitution into~\eqref{eq:Cn_Sn} proves~\eqref{eq:Cn_lagrange}.
\end{proof}

\end{document}